\documentclass{IEEEtran4PSCC}
\usepackage{algorithm, algpseudocode}
\usepackage{amsfonts, amsthm}
\usepackage{xcolor}
\usepackage{caption, subcaption}

\newtheorem{remark}{Remark}

\newcommand{\N}{\mathcal{N}} 
\newcommand{\E}{\mathcal{E}} 
\newcommand{\Et}{\mathcal{E}^{\tau}} 

\newcommand{\wc}{w^{\text{c}}}
\newcommand{\ws}{w^{\text{s}}}
\newcommand{\bc}{b^{\text{c}}}
\newcommand{\bs}{b^{\text{s}}}
\newcommand{\gs}{g^{\text{s}}}

\newcommand{\ffp}{f^{\text{p}}} 
\newcommand{\ffq}{f^{\text{q}}}
\newcommand{\UBgp}{\overline{f}^{\text{p}}} 
\newcommand{\UBgq}{\overline{f}^{\text{q}}} 
\newcommand{\LBgp}{\underline{f}^{\text{p}}} 
\newcommand{\LBgq}{\underline{f}^{\text{q}}}

\newcommand{\ddp}{d^{\text{p}}} 
\newcommand{\ddq}{d^{\text{q}}}

\newcommand{\lpp}{l^{\text{p+}}} 
\newcommand{\lpm}{l^{\text{p--}}} 
\newcommand{\lqp}{l^{\text{q+}}} 
\newcommand{\lqm}{l^{\text{q--}}} 
 
\newcommand{\Id}{I^{\text{d}}} 
\newcommand{\Ieff}{I^{\text{eff}}} 

\newcommand{\vd}{v^{\text{d}}} 

\newcommand{\Nd}{\mathcal{N}^{\text{d}}} 
\newcommand{\Ns}{\mathcal{N}^{\text{s}}} 
\newcommand{\Ed}{\mathcal{E}^{\text{d}}} 

\newcommand{\vr}{v^{\text{\tiny R}}}
\newcommand{\vi}{v^{\text{\tiny I}}}

\usepackage[noadjust]{cite}

\ifCLASSINFOpdf
   \usepackage[pdftex]{graphicx}
\else
   \usepackage[dvips]{graphicx}
\fi
\usepackage[cmex10]{amsmath}
\makeatletter
\let\old@ps@headings\ps@headings
\let\old@ps@IEEEtitlepagestyle\ps@IEEEtitlepagestyle
\def\psccfooter#1{%
    \def\ps@headings{%
        \old@ps@headings%
        \def\@oddfoot{\strut\hfill#1\hfill\strut}%
        \def\@evenfoot{\strut\hfill#1\hfill\strut}%
    }%
    \def\ps@IEEEtitlepagestyle{%
        \old@ps@IEEEtitlepagestyle%
        \def\@oddfoot{\strut\hfill#1\hfill\strut}%
        \def\@evenfoot{\strut\hfill#1\hfill\strut}%
    }%
    \ps@headings%
}
\makeatother

\psccfooter{%
        \parbox{\textwidth}{\hrulefill \\ \small{23rd Power Systems Computation Conference} \hfill \begin{minipage}{0.2\textwidth}\centering \vspace*{4pt} \includegraphics[scale=0.06]{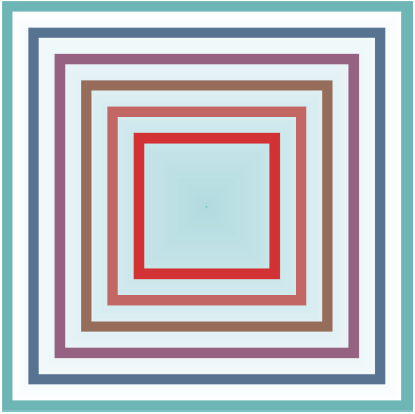}\\\small{PSCC 2024} \end{minipage} \hfill \small{Paris-Saclay, France --- June 4--7, 2024}}%
}

\usepackage[hyphens]{url}

\begin{document}
%
\title{Heuristic Algorithms for Placing Geomagnetically Induced Current Blocking Devices}


 


\author{\IEEEauthorblockN{Minseok Ryu\IEEEauthorrefmark{1},
Ahmed Attia\IEEEauthorrefmark{2},
Arthur Barnes\IEEEauthorrefmark{3}, 
Russell Bent\IEEEauthorrefmark{4}, 
Sven Leyffer\IEEEauthorrefmark{2}, and
Adam Mate\IEEEauthorrefmark{3}}

\vspace{3mm}
\IEEEauthorblockA{\IEEEauthorrefmark{1}School of Computing and Augmented Intelligence\\
Arizona State University, Tempe, AZ}
\vspace{1mm}
\IEEEauthorblockA{\IEEEauthorrefmark{2}Mathematics and Computer Science Division\\
Argonne National Laboratory, Lemont, IL}
\vspace{1mm}
\IEEEauthorblockA{\IEEEauthorrefmark{3}Information Systems and Modeling Group and \IEEEauthorrefmark{4}Applied Mathematics and Plasma Physics Group\\
Los Alamos National Laboratory, Los Alamos, NM}
}


\maketitle

\begin{abstract}
We propose a new heuristic approach for solving the challenge of determining optimal placements for geomagnetically induced current blocking devices on electrical grids. 
Traditionally, these determinations are approached by formulating the problem as mixed-integer nonlinear programming models and solving them using optimization solvers based on the spatial branch-and-bound algorithm.
However, computing an optimal solution using the solvers often demands substantial computational time due to their inability to leverage the inherent problem structure. 
Therefore, in this work we propose a new heuristic approach based on a three-block alternating direction method of multipliers algorithm, and we compare it with an existing stochastic learning algorithm. Both heuristics exploit the structure of the problem of interest. 
We test these heuristic approaches through extensive numerical experiments conducted on the EPRI-21 and UIUC-150 test systems. 
The outcomes showcase the superior performance of our methodologies in terms of both solution quality and computational speed when compared with conventional solvers. \\


\end{abstract}

\section{Introduction}
\begin{figure}[h!]
\centering
\includegraphics[width=.49\linewidth]{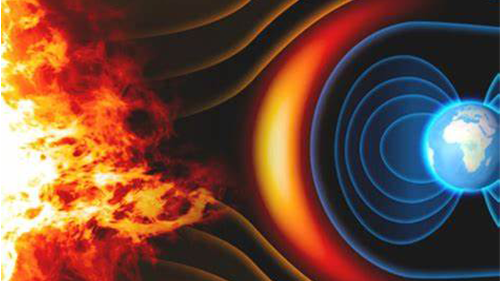}   
\includegraphics[width=.425\linewidth]{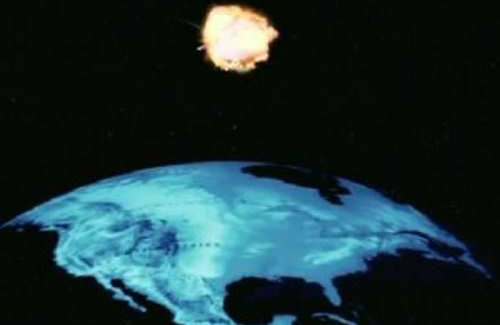}
\caption{Visualization of GMD (left) and  EMP (right) by the U.S. Department of Homeland Security \cite{dhs}.}
\label{fig:gmd_emp}
\end{figure}

This paper considers a geomagnetically induced current (GIC) blocking device placement problem, referred to herein as the GIC-BDP problem, for determining optimal locations for installing a limited number of devices to mitigate the adverse effect of GIC on transmission networks. 
GICs are low-frequency currents that can flow through transmission lines and transformers, which typically emerge because of a substantial electric field present on the Earth's surface, often referred to as the E-field.
The increase in E-field strength often stems from naturally occurring geomagnetic disturbances (GMDs) induced by severe space weather events or from intentional electromagnetic pulse (EMP) attacks, visualized in Fig.~\ref{fig:gmd_emp}.

The presence of GIC can give rise to various adverse effects, including the emergence of current harmonics, transformer saturation, and increased reactive power losses. Each of these has the potential to cause damage to critical equipment and even trigger cascading failures \cite{NERC2012-gmd, barnes21-hiddenfailures, mate21-pmsgmd-cascade}.
Considering the severity of these consequences, mitigating the flows of GIC holds significant importance in ensuring the resilience and reliability of bulk energy systems. A handful of potential strategies for GIC mitigation have been put forth in the literature: installation of blocking devices \cite{kappenman1991gic, zhu2014blocking, etemadi2014optimal, liang2015optimal, liang2019optimal} and employing transmission line switching \cite{lu2017optimal, ryu2020algorithms, ryu2022mitigating}.
The optimization of these mitigation approaches is often formulated as a mixed-integer nonlinear programming (MINLP) model. 
Optimal solutions to these models can be computed by using techniques such as the spatial branch-and-bound algorithm, as demonstrated in \cite{liang2015optimal}.  
This process can be time-consuming, however, mainly due to the presence of (i) binary variables that determine the placement of blocking devices or the selection of switched-on transmission lines, (ii) absolute-value equations used for computing the effective GIC, and (iii) nonlinear, nonconvex equations containing trilinear terms that involve trigonometric functions for computing AC optimal power flow (OPF) on transmission networks.

In this paper we reformulate the MINLP model for the GIC-BDP problem and propose a heuristic approach that provides good-quality solutions within reasonable computation time.   
Our contributions are summarized as follows:
\begin{enumerate}
\item We reformulate the absolute-value equations as complementarity constraints \cite{ferris1997engineering,fletcher2004solving,scheel2000mathematical} in the MINLP model to avoid potential convergence failure in MINLP solvers.  
\item We propose a three-block alternating direction method of multipliers (ADMM) algorithm for solving the MINLP model, where the first-block subproblem is an integer program that admits a closed-form solution and the second- and third-block subproblems are NLPs. 
\item For comparison, we adapt a stochastic learning (SL) algorithm for solving the MINLP model, which samples binary variables from a joint multivariate Bernoulli distribution whose parameters are optimized iteratively. 
\end{enumerate}
The remainder of this paper is organized as follows. 
In Section \ref{sec:model} we present our MINLP model for the GIC-BDP problem. 
In Section \ref{sec:heuristics} we present two heuristic approaches that exploit the structure of the problem.
In Section \ref{sec:experiments} we show that our heuristic approaches outperform the state-of-the-art MINLP solvers with respect to solution quality and computation speed.

\section{AC-OPF and GIC Blocker Model} \label{sec:model} 
Here we present a MINLP model for the GIC-BDP problem that determines optimal locations for installing a limited quantity of devices to effectively mitigate the detrimental effects of GIC on transmission networks.
In contrast to preceding studies \cite{lu2017optimal, ryu2022mitigating}, which relied on the polar representation of power flow, the rectangular form of power flow equation is embedded in our model to avoid potential numerical instability of MINLP solvers, primarily attributed to the presence of trilinear terms that involve trigonometric functions within the polar form. 
In what follows, we briefly describe each set of constraints in our MINLP model (see Table \ref{tab-definition} for notations, separated according to sets, parameters, and variables for AC-OPF and GIC).
\begin{table}[h!]
\footnotesize
\centering
\caption{Nomenclature}
\label{tab-definition}    
\begin{tabular}{l||l}
\hline
\multicolumn{2}{c}{\textbf{Sets and parameters}} \\ \hline
$\mathcal{G},\ \mathcal{N}, \ \mathcal{E}$ & set of generators, buses, and lines in ac network \\
$\mathcal{E}^{\tau} \subseteq \mathcal{E} $ & set of transformers \\
$\mathcal{E}_i, \mathcal{G}_i$ & set of lines and generators connected to $i \in \mathcal{N}$ \\		
$c^{\text{\tiny F1}}_k, c^{\text{\tiny F2}}_k$ & cost of generating power at $k \in \mathcal{G}$ \\
$\LBgp_k, \UBgp_k$ & bounds on the real power generation of $k \in \mathcal{G}$ \\
$\LBgq_k, \UBgq_k$ & bounds on the reactive power generation of $k \in \mathcal{G}$ \\
$\kappa$ & unit penalty cost for power unbalance at $i \in \mathcal{N}$ \\
$\ddp_i, \ddq_i$ & real and reactive power demand at $i \in \mathcal{N}$ \\
$\underline{\vr}_i, \overline{\vr}_i$, $\underline{\vi}_i, \overline{\vi}_i$ & voltage bounds at $i \in \mathcal{N}$ \\
$\gs_i, \bs_i$ & shunt conductance and susceptance at $i \in \mathcal{N}$ \\
$g_e, b_e$ & conductance, susceptance of $e \in \mathcal{E}$ \\
$\bc_e$ & line-charging susceptance of  $e \in \mathcal{E}$ \\		
$\overline{s}_e$ & apparent power limit of line $e \in \mathcal{E}$\\
$\underline{\theta}_{ij}, \overline{\theta}_{ij}$ & bounds on the phase angle difference at $e_{ij}  \in \mathcal{E}$ \\ 
$K_e$ & loss factor of transformer $e \in \mathcal{E}^{\tau}$\\
$\overline{I}_e$ & upper limit of the effective GIC on $e \in \mathcal{E}^{\tau}$\\
\hline 
$\Nd, \Ed$ & set of buses and lines in dc network \\
$\Ns \subseteq \Nd$ & set of substations \\
$\mathcal{E}^{\text{d}-}_m$ & set of incoming lines connected to $m \in \mathcal{N}^{\text{\tiny d}}$ \\ 
$\mathcal{E}^{\text{d}+}_m$ & set of outgoing lines connected to $m \in \mathcal{N}^{\text{\tiny d}}$ \\ 
$\gamma_{\ell}$ & conductance of $\ell \in \mathcal{E}^{\text{\tiny d}}$ \\
$a_m$ & inverse of ground resistance at $m \in \mathcal{N}^{\text{\tiny d}}$  \\		
$\xi_{\ell} $ & GIC-induced voltage sources on $\ell \in \mathcal{E}^{\text{\tiny d}}$ \\ 
\hline
\multicolumn{2}{c}{\textbf{Variables}} \\ \hline		
$\vr_i$, $\vi_i$ & real and imaginary part of complex voltage at $i \in \N$ \\
$\ffp_k, \ \ffq_k$ & real and reactive power generated by $k \in \mathcal{G}$\\
$p_{ei}, \ p_{ej}$ & real power flow on $e_{ij} \in \E$ at $i \in \N$ and $j \in \N$  \\
$q_{ei}, \ q_{ej}$ & reactive power flow on $e_{ij} \in \E$ at $i \in \N$ and $j \in \N$  \\		
$ \lpp_i, \lqp_i$ & real and reactive power load shedding at $i \in \N$ \\
$ \lpm_i, \lqm_i$ & real and reactive power overconsumed at $i \in \N$ \\		
$d^{\text{qloss}}_i$ & reactive power loss due to GIC at $i \in \N$ \\		
$\Ieff_e$ & effective GIC on $e \in \Et$ \\				
\hline 
$\Id_{\ell}$ & GIC that flows on $\ell \in \Ed$ \\
$\vd_m$ & GIC-induced voltage magnitude at $m \in \Nd$  \\		
$z_m \in \{0,1\}$ & $z_m=1$ if a device is installed at $m \in \Ns$ \\
\hline
\end{tabular}
\end{table}

\subsubsection{Operational constraints}
\begin{subequations}
\label{operational_constraints}  
\begin{align}
& p_{ei}^2 + q_{ei}^2 \leq  (\overline{s}_e)^2 , \ \ p_{ej}^2 + q_{ej}^2 \leq  (\overline{s}_e)^2, \ \forall e_{ij} \in \mathcal{E},  \label{thermal_limit} \\
& \ffp_k \in [\LBgp_k, \UBgp_k], \quad \ffq_k \in [\LBgq_k, \UBgq_k], \ \forall k \in \mathcal{G}, \label{bound_gen} \\ 
& \vr_i  \in [\underline{\vr}_i, \overline{\vr}_i], \ \ \vi_i  \in [\underline{\vi}_i, \overline{\vi}_i] \ \forall i \in \mathcal{N},  \label{bound_vol_mag}
\end{align}  
\end{subequations}
where Eq.~\eqref{thermal_limit}, Eq.~\eqref{bound_gen}, and Eq.~\eqref{bound_vol_mag} ensure line thermal limit, bounded power generation, and bounded voltage, respectively.

\subsubsection{Power flow equations}
For every line $e_{ij} \in \E$, we have
\begin{subequations}
\label{power_flow}
\begin{align} 
& p_{ei} =  g_e w_i - (g_e \wc_e + b_e \ws_e),  \\
& p_{ej} =  g_e w_j - (g_e \wc_e - b_e \ws_e),  \\
& q_{ei} = -(b_e + \bc_e/2) w_i + (b_e \wc_e - g_e \ws_e ), \\
& q_{ej} = -(b_e + \bc_e/2) w_j + (b_e \wc_e + g_e \ws_e ), 
\end{align}
where
\begin{align}
& w_i = \vr_i \vr_i + \vi_i \vi_i , \ \forall i \in \N, \\
& \wc_e = \vr_i \vr_j + \vi_i \vi_j, \ \forall e_{ij} \in \E, \\
& \ws_e = \vr_j \vi_i - \vr_i \vi_j, \ \forall e_{ij} \in \E, \\
& \tan (\underline{\theta}_{ij}) \wc_e \leq \ws_e \leq \tan (\overline{\theta}_{ij}) \wc_e, \label{angle_tangent}
\end{align}
\end{subequations}
which ensure that power flow is governed by Ohm's law.  

\subsubsection{Balance equations with reactive power losses by GIC}
For every bus $i \in \N$, we have
\begin{subequations}
\label{balance}
\begin{align}
& \sum_{e \in \mathcal{E}_i } p_{ei} = \sum_{k \in \mathcal{G}_i } \ffp_k - \ddp_i + \lpp_i - \lpm_i - \gs_i w_i,  \\
& \sum_{e \in \mathcal{E}_i } q_{ei} = \sum_{k \in \mathcal{G}_i } \ffq_k - \ddq_i  + \lqp_i - \lqm_i + \bs_i w_i - d^{\text{qloss}}_i,  
\end{align}
\end{subequations}
which ensure power balance.
Note that $d^{\text{qloss}}_i$ represents reactive power losses by GIC, computed by
\begin{subequations}
\label{dqloss_polar}
\begin{align}
& d^{\text{qloss}}_i =  \sum_{e \in \mathcal{E}^{\tau}_i} K_e \sqrt{v_i} \Ieff_e, \ \forall i \in \mathcal{N}, \label{dqloss_polar_1} \\
& \Ieff_e \in [0, \overline{I}_e], \ \forall e \in \Et,
\end{align}
\end{subequations}
where $\Ieff_e$ is the amounts of effective GIC at a transformer $e \in \Et$, which are computed utilizing a representative DC network ($\Nd, \Ed$) derived from modifications made to the original AC network ($\N, \E$), as shown in Fig.~\ref{fig:ACDC}.

\subsubsection{GIC Model}
\begin{figure}[!h]
  \centering
  \includegraphics[scale=0.4]{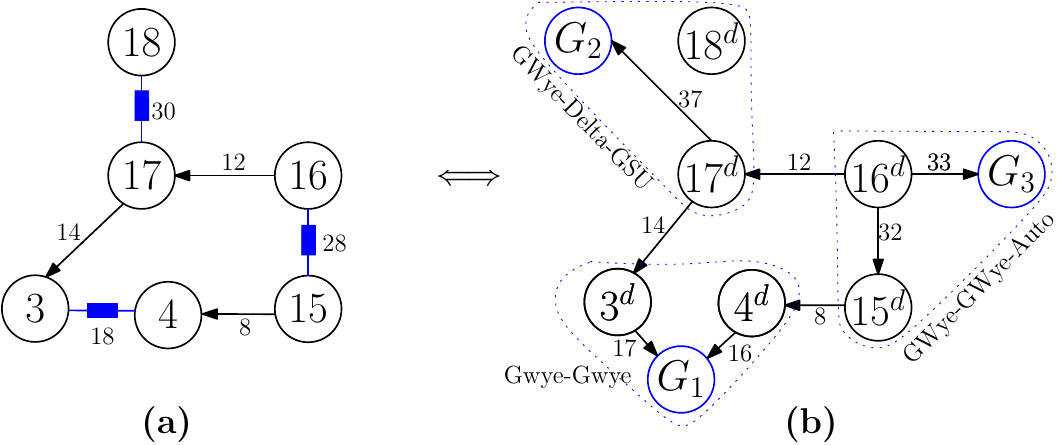}
  \caption{AC (left) and DC (right) power network \cite{ryu2022mitigating}. }
  \label{fig:ACDC} 
\end{figure}

The DC network is constructed by adding a set of substations (e.g., $G_1$, $G_2$, and $G_3$ in Fig.~\ref{fig:ACDC}) to the underlying AC network.
For each line of the DC network, the GIC-induced voltage source is given by $\xi_{\ell}$, which is zero for all transformers (e.g., line numbers 16, 17, 32, 33, and 37 in Fig.~\ref{fig:ACDC}) and has nonzero values for all transmission lines. With this information, the GIC can be computed as
\begin{align}
& \Id_{\ell} = \gamma_{\ell} (\vd_m - \vd_n + \xi_{\ell}), \ \forall \ell_{mn} \in \Ed, \label{gic} 
\end{align}
where $\vd$ is the GIC-induced voltage magnitude.

Installing a blocking device disconnects the transformer neutral from a substation and changes the conductance matrix. To model this, the GIC balance equations are introduced:
\begin{align}
& \sum_{\ell \in \mathcal{E}^{d-}_m} \Id_{\ell} - \sum_{\ell \in \mathcal{E}^{d+}_m} \Id_{\ell} = a_m \vd_m (1-z_m), \ \forall m \in \Nd, \label{gic_balance}
\end{align} 
which deduct the effect of $a_m$ if a device is installed at $m$.
For notation brevity, we set $a_m = 0$ for $m \in \Nd \setminus \Ns$ to ensure that the binary variables $z_m$ are  defined only for $m \in \Ns$. 
To limit the number of blocking devices, we introduce a bound: 
\begin{subequations}
\label{budget_constraint} 
\begin{align}
& \sum_{m \in \Ns} z_m \leq V,  \\
& z_m \in \{0,1\}, \ \forall m \in \Ns.
\end{align}
\end{subequations} 

\subsubsection{Effective GIC computation}
The GIC in Eq.~\eqref{gic} is used to calculate the effective GIC in Eq.~\eqref{dqloss_polar_1} through the following absolute-value equations for different types of transformers:
\begin{align}
& \Ieff_e = | \Theta_e |, \ \forall e \in \Et, \label{eff_gic_abs}
\end{align}
where 
\begin{subequations}
\label{Theta}  
\begin{align}
& \Theta_e =  \big(\frac{N_h \Id_h + N_l \Id_l}{N_h} \big),   \text{ if $e$ is GWye-GWye}, \label{Gwye} \\
& \Theta_e =   \big(\frac{N_s \Id_s + N_c \Id_c}{N_s + N_c} \big),   \text{ if $e$ is GWye-GWye Auto} , \label{Auto}\\
& \Theta_e = \Id_h,   \text{ if $e$ is GWye-Delta-GSU}, \label{GSU} 
\end{align}
\end{subequations}
where $N_h$, $N_l$, $N_s$, and $N_c$ are the number of turns in the high-side, low-side, series, and common windings, respectively. 

The absolute value in Eq.~\eqref{eff_gic_abs} is nonsmooth, which can cause convergence failure in MINLP solvers. 
To avoid this situation, we apply the equivalent smooth complementarity reformulation \cite{fletcher2004solving, ferris1997engineering,scheel2000mathematical, raghunathan2005interior}:
\begin{subequations}
\label{eff_gic_smooth}  
\begin{align}
    & s_e^+ - s_e^- = \Theta_e, \ \forall e \in \Et, \\
    & \Ieff_e = s_e^+ + s_e^-,  \ \forall e \in \Et, \label{eff_gic_smooth_2}  \\
    & s_e^- \geq 0, \ \ s_e^+ \geq 0, \ \  s_e^+ s_e^- \leq 0, \ \forall e \in \Et, 
\end{align}
\end{subequations}
where $s_e^+$ and $s_e^-$ are slack variables.
 
\subsubsection{MINLP models}
Now we can define our MINLP model: 
\begin{align}
& \hspace{-8mm} \texttt{(AC-rect)} \nonumber \\
\min \ & \sum_{k \in \mathcal{G}} ( c^{\text{F1}}_k \ffp_k +  c^{\text{F2}}_k (\ffp_k)^2 )  + \sum_{i \in \mathcal{N}} \kappa ( \lpp_i + \lpm_i + \lqp_i + \lqm_i  )  \nonumber \\
\mbox{s.t.} \ 
& Eq.~\eqref{operational_constraints} - Eq.~\eqref{budget_constraint}, 
Eq.~\eqref{Theta}, 
Eq.~\eqref{eff_gic_smooth}. 
\nonumber 
\end{align}


\section{Heuristic Algorithms} \label{sec:heuristics}

In this section we propose a new heuristic approach for the MINLP model described in Section~\ref{sec:model}, which we call a three-block ADMM algorithm with a binary subproblem, which we will refer to as 3ADMM-B. We also show how a stochastic learning  approach can be applied to our problem for comparison purposes.

\subsection{Three-Block Alternating Direction Method of Multipliers}
ADMM \cite{douglas1956numerical,eckstein1992douglas} is typically applied to convex, continuous optimization problems that can be decomposed into two blocks. The method  is related to augmented Lagrangian methods and consists of solving a sequence of alternating optimization problems followed by a first-order multiplier update. In contrast, we apply ADMM to a discrete optimization problem.
The proposed 3ADMM-B is derived by exploiting the structure of the problem; thus it is a problem-specific algorithm. 

First, we observe that the AC and DC network formulations are connected through the effective GIC variables $\{\Ieff_e\}_{e \in \Et}$ in Eq.~\eqref{dqloss_polar} and Eq.~\eqref{eff_gic_smooth_2}. 
By introducing auxiliary variables $I^{\text{ac}}$ and $I^{\text{dc}}$, we separate  constraints for the AC network from those for the DC network:
\begin{subequations}
\begin{align}
& g^{\text{ac}} (x, I^{\text{ac}}) \leq 0, \label{admm_form_ac} \\
& g^{\text{dc}} (y, z, I^{\text{dc}}) \leq 0, \label{admm_form_dc}  \\  
& I^{\text{dc}}_e = I^{\text{ac}}_e \in [0, \overline{I}_e], \ \forall e \in \mathcal{E}^{\tau}, \label{consensus}  
\end{align}
\end{subequations}
where Eq.~\eqref{admm_form_ac} and Eq.~\eqref{admm_form_dc} represent constraints for AC and DC networks, respectively, and Eq.~\eqref{consensus} represents consensus constraints.
Note that $x$ and $y$ are continuous local variables while $z$ is a binary vector that should satisfy Eq.~\eqref{budget_constraint}.

Second, we introduce auxiliary variables $z^{\text{b}}$ to remove the binary restriction from the DC network:
\begin{subequations}
\begin{align}
& z^{\text{b}}_i = z_i, \ \forall i \in [S], \label{binary_consensus} \\
& z^{\text{b}}_i  \in \{0,1\}, \ \ z_i \in [0,1], \ \forall i \in [S], \label{binary_aux} \\
& \sum_{i=1}^S z^{\text{b}}_i \leq V, \label{binary_budget}
\end{align}
\end{subequations}
where $S := |\Ns|$.
We note that the consensus constraints Eq.~\eqref{binary_consensus} ensure that the continuous copy agrees with the binary choices.
The MINLP model is written in the following form:
\begin{align}
\min \ & f(x) \nonumber \\
\mbox{s.t.} \  
& Eq.~\eqref{admm_form_ac} - Eq.~\eqref{consensus}, Eq.~\eqref{binary_consensus} - Eq.~\eqref{binary_budget}, \nonumber
\end{align}
where $f$ corresponds to the objective function of 
\texttt{(AC-rect)}.

By introducing dual variables $\lambda$ and $\mu$ associated with Eq.~\eqref{binary_consensus} and Eq.~\eqref{consensus}, respectively, the augmented Lagrangian is given by 
\begin{align}
\max_{\lambda, \mu} \ \min_{z,z^b,I} \ & f(x) + \langle \lambda, z^{\text{b}} - z \rangle + \langle \mu, I^{\text{dc}}- I^{\text{ac}} \rangle \label{augLag} \\
& + \frac{\rho}{2} \Big\{ \| z^{\text{b}} - z  \|^2 + \| I^{\text{dc}}- I^{\text{ac}}  \|^2 \Big\} \nonumber \\
& \text{s.t.}\; Eq.~\eqref{admm_form_ac}, Eq.~\eqref{admm_form_dc}, Eq.~\eqref{binary_aux}, Eq.~\eqref{binary_budget} , \nonumber
\end{align}
and we now consider an approach that solves \eqref{augLag} instead of \texttt{(AC-rect)}.
In the $t$th iteration of the proposed 3ADMM-B, the first-block subproblem is given by
\begin{align}
\min_{ z^{\text{b}} \in \{0,1\}^S} \ & \langle \lambda^{(t)}, z^{\text{b}} \rangle + \frac{\rho}{2} \| z^{\text{b}} - z^{(t)}\|^2  \label{admm_first_block} \\
\text{s.t.}\quad  & Eq.~\eqref{binary_budget} , \nonumber
\end{align}
the second-block subproblem is given by 
\begin{align}
\min_{z \in [0,1]^S, I^{\text{dc}} \in [0,\overline{I}]} \ & - \langle \lambda^{(t)}, z \rangle + \langle \mu^{(t)}, I^{\text{dc}} \rangle  \label{admm_second_block}  \\
&+ \frac{\rho}{2} \Big\{ \| z^{\text{b}(t+1)} - z \|^2+ \| I^{\text{dc}} - I^{\text{ac}(t)} \|^2 \Big\} \nonumber \\
\text{s.t.}\quad & Eq.~\eqref{admm_form_dc}, \nonumber 
\end{align}
the third-block subproblem is given by 
\begin{align}
\min_{ I^{\text{ac}} \in [0,\overline{I}]} \ & f(x) - \langle \mu^{(t)}, I^{\text{ac}} \rangle  + \frac{\rho}{2} \| I^{\text{dc}(t+1)} - I^{\text{ac}} \|^2 \label{admm_third_block}  \\
\text{s.t.} \quad & Eq.~\eqref{admm_form_ac}, \nonumber
\end{align}
and the dual update is given by
\begin{subequations}
\label{admm_dual_update}  
\begin{align}
& \lambda^{(t+1)} = \lambda^{(t)} + \rho(z^{\text{b}(t+1)} - z^{(t+1)}), \\
& \mu^{(t+1)} = \mu^{(t)} + \rho(I^{\text{dc}(t+1)} - I^{\text{ac}(t+1)}).
\end{align}
\end{subequations}
We note that Eq.~\eqref{admm_second_block} and Eq.~\eqref{admm_third_block} are NLP models with convex quadratic objective functions that are easy to solve, while Eq.~\eqref{admm_first_block} is a convex quadratic program with binary variables that is easy to solve, as pointed out in the next remark.

\begin{remark}
Since $z$ is binary, we have $ \|z \|^2 = \langle \textbf{1}, z \rangle$, where $\textbf{1}$ is a vector with all components being $1$.
Therefore, Eq.~\eqref{admm_first_block} can be rewritten as follow:  
\begin{align}
\min_{z \in \{0,1\}^S} \ & \sum_{i=1}^{n} \big( \frac{\rho}{2} + \lambda^{(t)}_i - \rho z^{(t)}_i \big) z_i \nonumber \\
\text{s.t.} \quad & Eq.~\eqref{binary_budget}, \nonumber
\end{align}
which is a binary knapsack problem whose constraint coefficients are all one. Therefore, an optimal solution can be obtained greedily, as described in lines 9--17 of Alg.~\ref{alg:admm}.
Specifically, we first sort the objective coefficient in an increasing order, namely, $\widehat{c}_{i_1} \leq \ldots \leq \widehat{c}_{i_{S}}$, and sequentially set $z
_j=1$ if $\widehat{c}_{i_j} < 0$ for $j \in \{i_1, \ldots, i_S\}$ until the budget $V$ is consumed.
\end{remark}

In Alg.~\ref{alg:admm} we describe the proposed 3ADMM-B algorithm, composed of the three subproblems and dual updates as in lines 3, 4, 5, and 6, respectively.
We set the termination criterion based on the normalized primal and dual residuals, as described in \cite{wohlberg2017admm}. That is, the algorithm is terminated at iteration $t$ if the following holds:
\begin{align}
\max( p^{(t)}, d^{(t)} ) < \epsilon,
\end{align}
where $\epsilon$ is the tolerance of our algorithm and $p^{(t)}$ and $d^{(t)}$ are normalized primal and dual residuals computed at the $t$th iteration, respectively, given as
\begin{align*}
& p^{(t)} := \frac{\|v^{(t)}-u^{(t)}\|}{\max( \|u^{(t)}\|, \|v^{(t)}\|) }, \ \ d^{(t)} :=\frac{\rho\|u^{(t)}-u^{(t-1)}\|}{  \|w^{(t)}\| } \\
& v^{(t)} := 
\begin{bmatrix}
z^{\text{b}(t)} \\ I^{\text{dc}(t)}
\end{bmatrix}, \ \ 
u^{(t)} := 
\begin{bmatrix}
z^{(t)} \\ I^{\text{ac}(t)}
\end{bmatrix}, \ \
w^{(t)} := 
\begin{bmatrix}
\lambda^{(t)} \\ \mu^{(t)}
\end{bmatrix}.
\end{align*}
  
\begin{algorithm}
\caption{Three-block ADMM with bianry (3ADMM-B)}\label{alg:admm}
\begin{algorithmic}[1]
\State Initialization: $t \gets 0$, $\lambda^{(t)}$, $\mu^{(t)}$, $z^{(t)}$, $I^{\text{ac}(t)}$
\While{not converged}
\State Compute $z^{\text{b}(t+1)} \gets \texttt{closed}(\rho, \lambda^{(t)}, z^{(t)})$ 
\State Compute $z^{(t+1)}$ and $I^{\text{dc}(t+1)}$ by solving Eq.~\eqref{admm_second_block}
\State Compute $I^{\text{ac}(t+1)}$ by solving Eq.~\eqref{admm_third_block}
\State Update duals by Eq.~\eqref{admm_dual_update}
\EndWhile
\State Return $z^{\text{b}(T+1)}$
\Statex \hrulefill
\Statex \texttt{closed}($\rho, \lambda, z^{c}$): 
\State Initialization: $z_i = 0$ for all $i \in [S]$ and budget $B \gets V$
\State Define $\widehat{c}_i := \frac{\rho}{2} + \lambda_i - \rho z^{\text{c}}_i$ for all $i \in [S]$
\State Sort elements of $\widehat{c}$ such that $\widehat{c}_{i_1} \leq \ldots \leq \widehat{c}_{i_{S}}$
\For{$j \in \{i_1, \ldots, i_S\}$}
\If{$\widehat{c}_{j} < 0 \text{ and } B > 0 $}
\State $z_j \gets 1$ and $B \gets B-1$
\EndIf
\EndFor
\State Return $z$
\end{algorithmic} 
\end{algorithm}

\subsection{Stochastic learning approach}

We also apply stochastic learning for binary optimization \cite{attia2022stochastic} to our problem. This heuristic finds binary solutions by sampling from a joint multivariate Bernoulli distribution whose probabilities are updated iteratively. It has the advantage that we can easily sample from the final distribution to explore possible alternative solutions, whereas 3ADMM-B is deterministic and produces only a single solution.

To describe the SL approach, we write \texttt{(AC-rect)} as
\begin{subequations}
\label{SL_Model}
\begin{align}
\min_{z \in \mathcal{Z}}  \ & F(z) ,
\end{align}
where
\begin{align}
& \mathcal{Z} := \Big\{ z \in \{0,1\}^S : \sum_{i=1}^S z_i \leq V \Big\}
\end{align}
\end{subequations}
and $F(z)$ is the optimal value of an NLP model resulting from fixing binary variables in \texttt{(AC-rect)} to some value $z \in \mathcal{Z}$.

The existing SL approach \cite{attia2022stochastic} has been developed for solving $\min_{z \in \{0,1\}^S} F(z)$ (e.g., Eq.~\eqref{SL_Model} without the constraint): 
\begin{align}
p^* = \arg \min_{p \in [0, 1]^S } \Phi(p) & := \mathbb{E}_{z \sim \mathbb{P}(z|p)} [F(z)] \nonumber \\
& = \sum_{k=1}^{2^S} P(\hat{z}^k|p) F(\hat{z}^k), \label{sl_model}
\end{align}
where $\mathbb{P}(z|p)$ represents a joint multivariate Bernoulli distribution with the probability mass function, $P(z|p) := \prod_{i=1}^{S}  p_i^{z_i} (1-p_i)^{1-z_i}$.
Equation~\eqref{sl_model} is different from the original problem  in that it aims to optimize probabilities $\{p_1, \ldots, p_S\}$ associated with binaries $\{z_1, \ldots, z_S \}$, leading to heuristic solutions.
Also, Eq.~\eqref{sl_model} can be considered as a machine learning model with $2^S$ number of data points. 
Thus, one can utilize stochastic gradient descent (SGD) types of algorithms for solving Eq.~\eqref{sl_model}.
 For more details on the existing SL approach, we refer the reader to \cite{attia2022stochastic}.

In this application, however, the existing approach cannot be immediately utilized because $z$ sampled from the distribution may not satisfy the budget constraint $\sum_{i=1}^S z_i \leq V$.
To address this issue, we first sort the probabilities in a decreasing order, namely, $p_{i_1} \geq \ldots \geq p_{i_S}$, and sample $z_j$ for all $j \in \{i_1, \ldots, i_S\}$ until the budget $V$ is consumed. 
The proposed sampling is described in lines 3 and 9--18 of Alg.~\ref{alg:stochastic}.
We set the termination criterion based on the norm of gradient, namely, $g^{(t)} < \epsilon$, where $\epsilon$ is a tolerance level; and we set the diminishing step size $\eta^{(t)}=a/t$, where $a > 0$ is some constant.

\begin{algorithm}[h!]
\caption{Stochastic learning approach}\label{alg:stochastic}
\begin{algorithmic}[1]
\State Initialization: probability $p^{(1)}$,
step size $\eta^{(1)}$, sample size $N$, and set $t \gets 1$
\While{not converged}
\State Sample $N$ scenarios of $\hat{z} \gets \texttt{sample} (p^{(t)})$
\State Compute a gradient:
\begin{align}
g^{(t)} \gets \frac{1}{N} \sum_{k=1}^N F (\hat{z}^k) \Big\{ \sum_{i=1}^{S} \Big( \frac{\hat{z}^k_i}{p^{(t)}_i} - \frac{1 - \hat{z}^k_i}{1-p^{(t)}_i} \Big) \vec{e}_i \Big\} \label{gradient_sl}
\end{align}
\State Update
\begin{align}
& p^{(t+1)} \gets \text{Proj}_{[0,1]^{S}} ( p^{(t)} - \eta^{(t)} g^{(t)} )  \\
& t \gets t+1 \nonumber
\end{align}
\EndWhile
\State Sample $N$ scenarios of $\hat{z} \gets \texttt{sample} (p^{(t+1)})$
\Statex \hrulefill
\Statex \texttt{Sample}($p$): 
\State Initialization: $z_i = 0$ for all $i \in [S]$ and budget $B \gets V$.
\State Sort elements of $p$ such that $p_{i_1} \geq \ldots \geq p_{i_S}$
\For{$j \in \{i_1, \ldots, i_S\}$}
\State Sample $z_j$ from the Bernoulli distribution $\mathbb{P}(z|p_{j})$
\If{$z_j = 1$}
\State $B \gets B-1$
\EndIf
\State Break if $B = 0$
\EndFor
\State Return $z$
\end{algorithmic}  
\end{algorithm}

\section{Numerical Experiments} \label{sec:experiments}

In this section we numerically demonstrate that the proposed heuristic approaches provide good-quality solutions within reasonable computation time, much faster than existing MINLP solvers (e.g., SCIP \cite{BestuzhevaEtal2021ZR} and Juniper \cite{kroger2018juniper}).
To achieve this, we employ our proposed approach and established solvers to compute solutions for \texttt{(AC-rect)} within a maximum time limit of 1 hour. We then assess solution quality by solving an NLP model that arises from fixing binary variables in the MINLP model to the obtained solutions.
Based on the EPRI-21 and UIUC-150 test systems, as described in \cite{ryu2022mitigating}, we designed case studies by varying the magnitude of E-field $E \in \{5, 10, 15, 20\}$ V/km while keeping the E-field direction fixed at 45 degrees.
For the EPRI-21 and UIUC-150 system, which have $8$ and $98$ substations, respectively, each of these substations is a potential location for installing a blocking device.
For illustration of the methodology, we set the budget as $V=3$ for EPRI-21 and $V=30$ for UIUC-150, with the choice of $V$ proportional to the size of the network.

We used Julia 1.8.1. for writing \texttt{(AC-rect)} via JuMP \cite{Lubin2023} and implementing the heuristic approaches.
For all experiments, we used Bebop, a 1024-node computing cluster (each computing node has 36 cores with Intel Xeon E5-2695v4 processors and 128 GB DDR4 of memory) at Argonne National Laboratory.
With this computing resource, we solved \texttt{(AC-rect)} using our heuristic approach and well-established solvers, each of which is executed sequentially.

\subsection{Motivating examples}

To motivate the development of heuristic approaches, we demonstrate that the optimal placement of blocking devices allows the power system to operate during severe GMDs with a smaller number of blockers, significantly reducing placement cost. We also show that this problem is computationally challenging mainly because of the presence of binary variables.

\begin{figure}[h!]
\centering      
\includegraphics[width=0.49\linewidth]{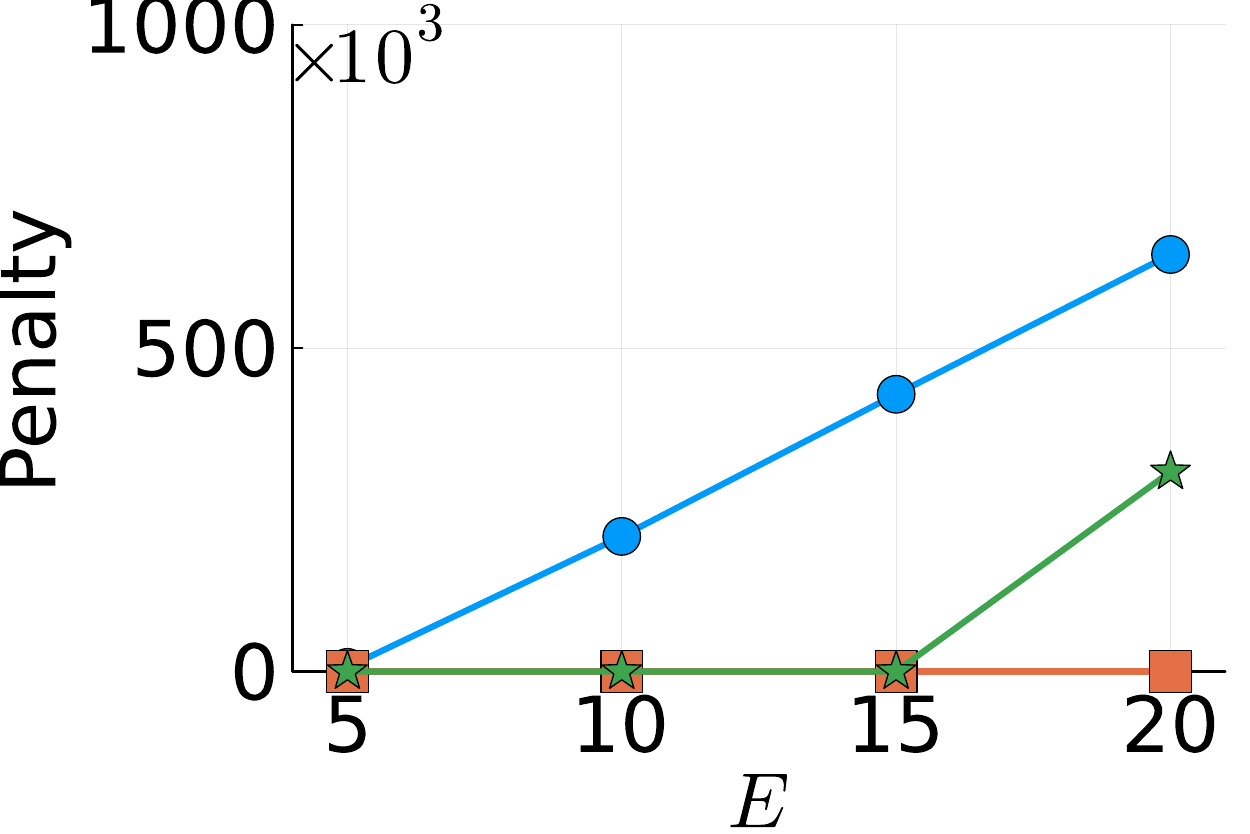} 
\includegraphics[width=0.49\linewidth]{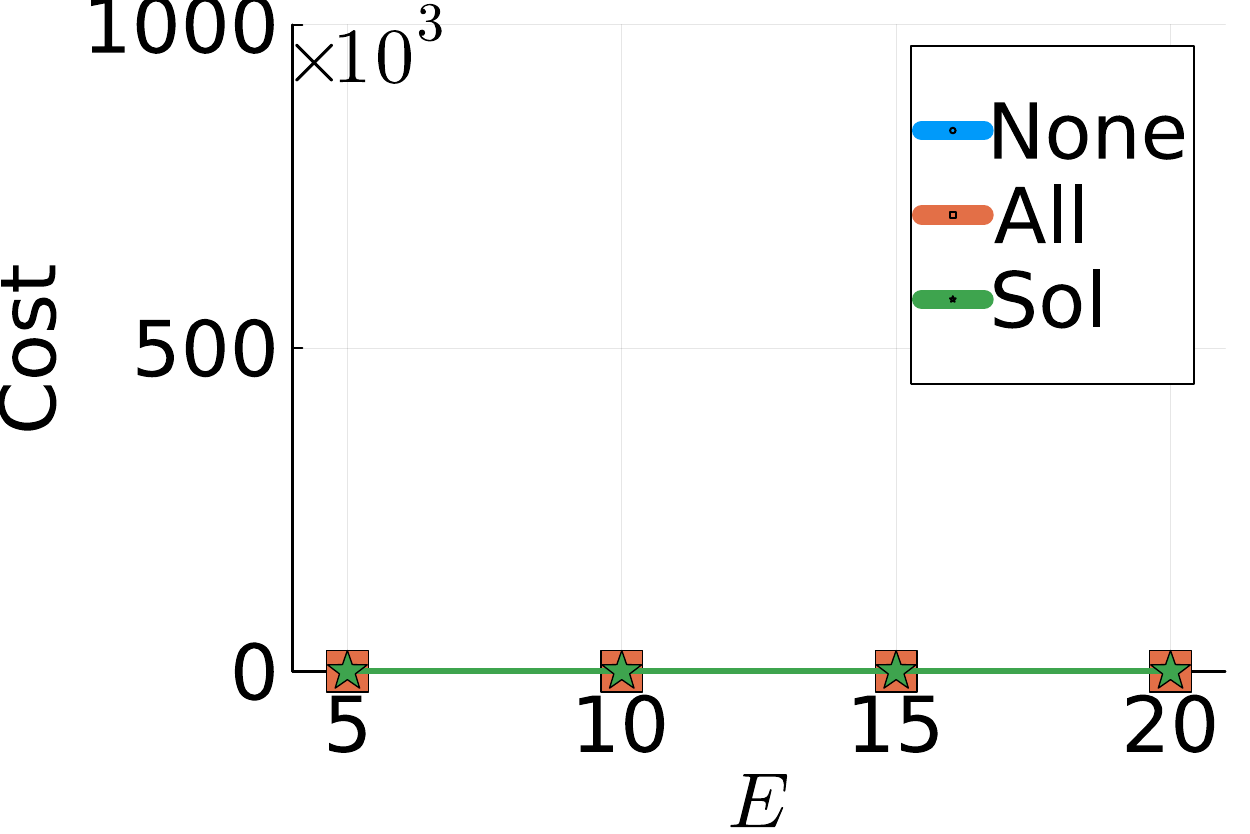} 
\caption{Effect of installing devices on the load-shedding penalty (left) and the power generation cost (right).}
\label{fig:effect_of_blockers}
\end{figure}

First, we employ the open-source MINLP solver SCIP to solve \texttt{(AC-rect)} using instances constructed by varying the E-field magnitude $E$ within the EPRI-21 system.
The results, consisting of the selected substations for installing blocking devices, are reported in Table~\ref{tab:solutions_epri21_scip}.
By solving NLP models derived from fixing the binary variables in \texttt{(AC-rect)} to the obtained solutions, we compute the load-shedding penalty and power generation cost in Fig.~\ref{fig:effect_of_blockers} (labeled as ``Sol'').
This is compared with scenarios where no blocking device is present (``None''), as well as when all substations have blocking devices (``All'').
For each $E$ value, installing devices on substations indicated in Table~\ref{tab:solutions_epri21_scip} mitigates the load-shedding penalty while preserving the power generation cost. This signifies that the solutions produced by SCIP within the 1-hour time limit can effectively alleviate the adverse effects of GIC on the power grid, despite not being optimal. 
We note that although installing blocking devices for all substations eliminates load shedding, the associated installation expenses can be significant. This situation highlights the necessity of identifying optimal sites for placing these devices.

\begin{figure}[h!]
\centering      
\includegraphics[width=.49\linewidth]{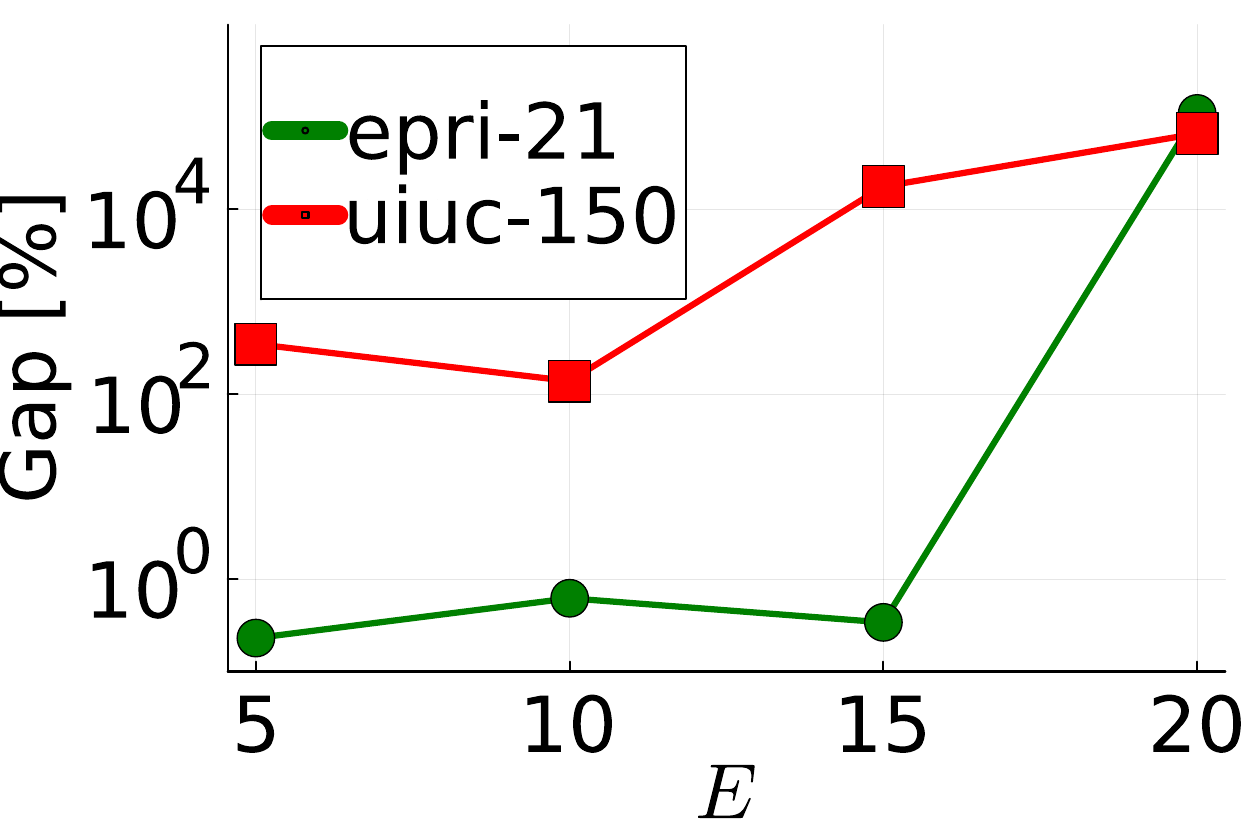} 
\caption{Gap computed at the time limit of 1 hour. }
\label{fig:gap}
\end{figure}

Unfortunately, solving these MINLP models to optimality is computationally intractable.
The solutions reported in Table~\ref{tab:solutions_epri21_scip} are incumbent solutions only obtained within the imposed time limit.
The corresponding solution gaps are outlined in Fig.~\ref{fig:gap}.
Notably, for the EPRI-21 system, smaller gaps are apparent when $E \leq 15$, while the gap becomes considerably larger at $E = 20$. For the larger UIUC-150 system, the gap is even more pronounced. This underscores the necessity to develop heuristic methods capable of generating high-quality solutions within reasonable computational time.

\begin{table}[h!]
\centering
\caption{Set of substations to install blocking devices produced by solving \texttt{(AC-rect)} using SCIP.}
\begin{tabular}{c|c|c|c|c}
$E$ [V/km] & $5$ & $10$ & $15$ & $20$ \\
\hline 
Substations & $\{3, 8\}$ & $\{5, 6, 8\}$& $\{1, 6, 8\}$ & $\{2, 6\}$ \\
\end{tabular}
\label{tab:solutions_epri21_scip}
\end{table}
 
\subsection{Comparison of heuristic approaches}

\begin{figure}[h!]
\centering        
\includegraphics[width=.49\linewidth]{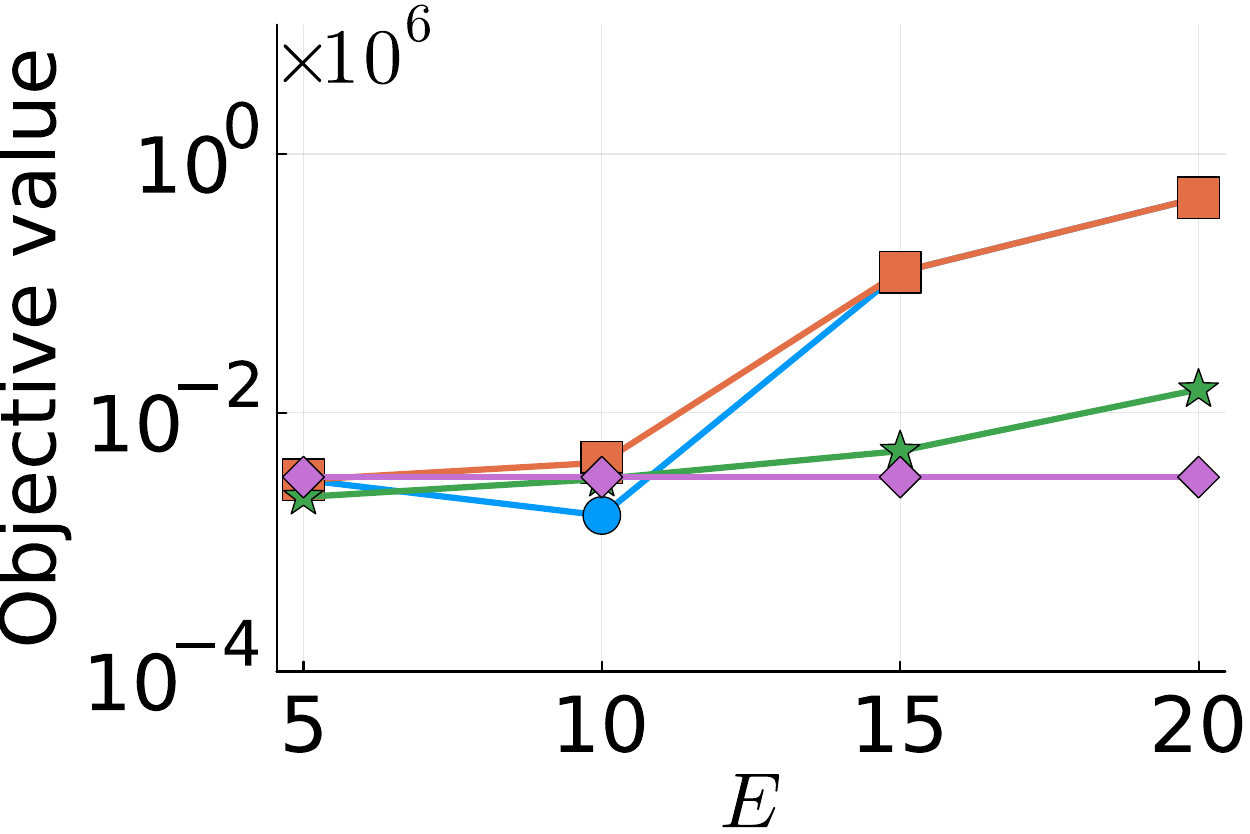}
\includegraphics[width=.49\linewidth]{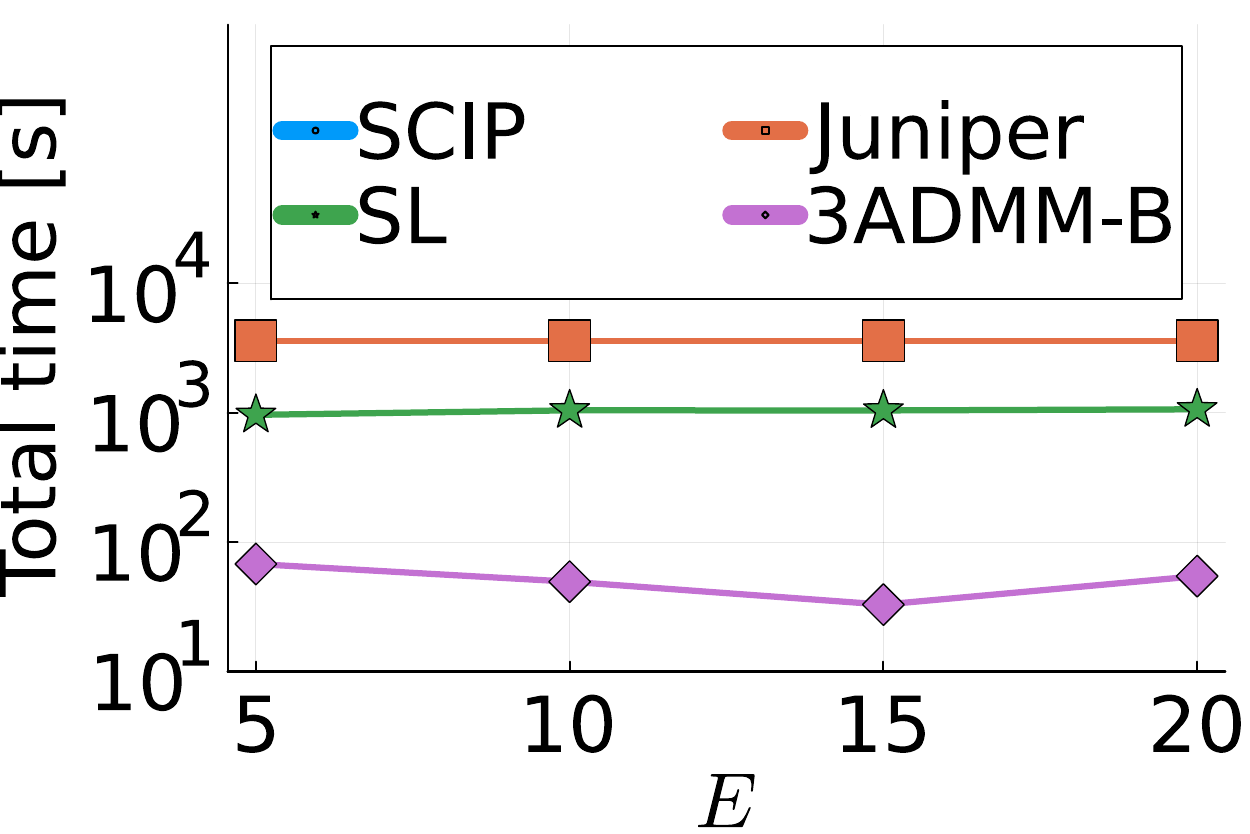}
\caption{Comparison of heuristic approaches with respect to the objective value (left) and computation time (right). }
\label{fig:uiuc150}
\end{figure}

In this section we numerically show that the proposed heuristic approaches provide superior solution quality compared with existing MINLP solvers.
To this end, we first solve \texttt{(AC-rect)} constructed by varying $E \in \{5, 10, 15, 20\}$ on the UIUC-150 system by our approaches (i.e.,  3ADMM-B and SL), SCIP, and Juniper, respectively.
To compare quality of solutions (i.e., where to place devices), we solve an NLP model derived from fixing the binary variables in \texttt{(AC-rect)} using the solutions, and we then report the resulting objective values in Fig.~\ref{fig:uiuc150}.
Additionally, we provide information on the total time required for computation and evaluation in the same figure.
The outcomes reveal that the solutions obtained through our heuristic approaches exhibit lower objective values in comparison with solutions generated by SCIP and Juniper, all within a 1-hour time limit. Notably, the enhanced quality of solutions by 3ADMM-B and SL is achieved within computation times of less than 100 and 1000 seconds, respectively. While the 3ADMM-B approach is generally faster than stochastic learning, the stochastic learning approach provides benefits in situations where it is necessary to explore the solution space of near-optimal solutions.
Further details regarding these two approaches will be discussed in the subsequent sections.

\subsection{Details on 3ADMM-B}

In this section we provide details on 3ADMM-B.
The ADMM penalty parameter $\rho$ is a hyperparameter that should be tuned for better performance in practice. 
We utilize the normalized residual balancing (NRB) technique \cite{wohlberg2017admm} to adaptively choose the value of $\rho$ in every iteration $t$ of the algorithm. 
Specifically, for given $\beta, \tau \in \mathbb{R}_+$ and the primal and dual residuals $p^{(t)}, d^{(t)}$, we update $\rho$ as follows:
\begin{align}
\rho^{(t+1)} \gets
\begin{cases}
\rho^{(t)} \tau    & \ \ \text{ if } p^{(t)} > \beta d^{(t)}     \\
\rho^{(t)} / \tau    & \ \ \text{ if } p^{(t)} < \beta d^{(t)}.    
\end{cases}
\end{align}

\begin{figure}[h!]
\centering      
\includegraphics[width=.49\linewidth]{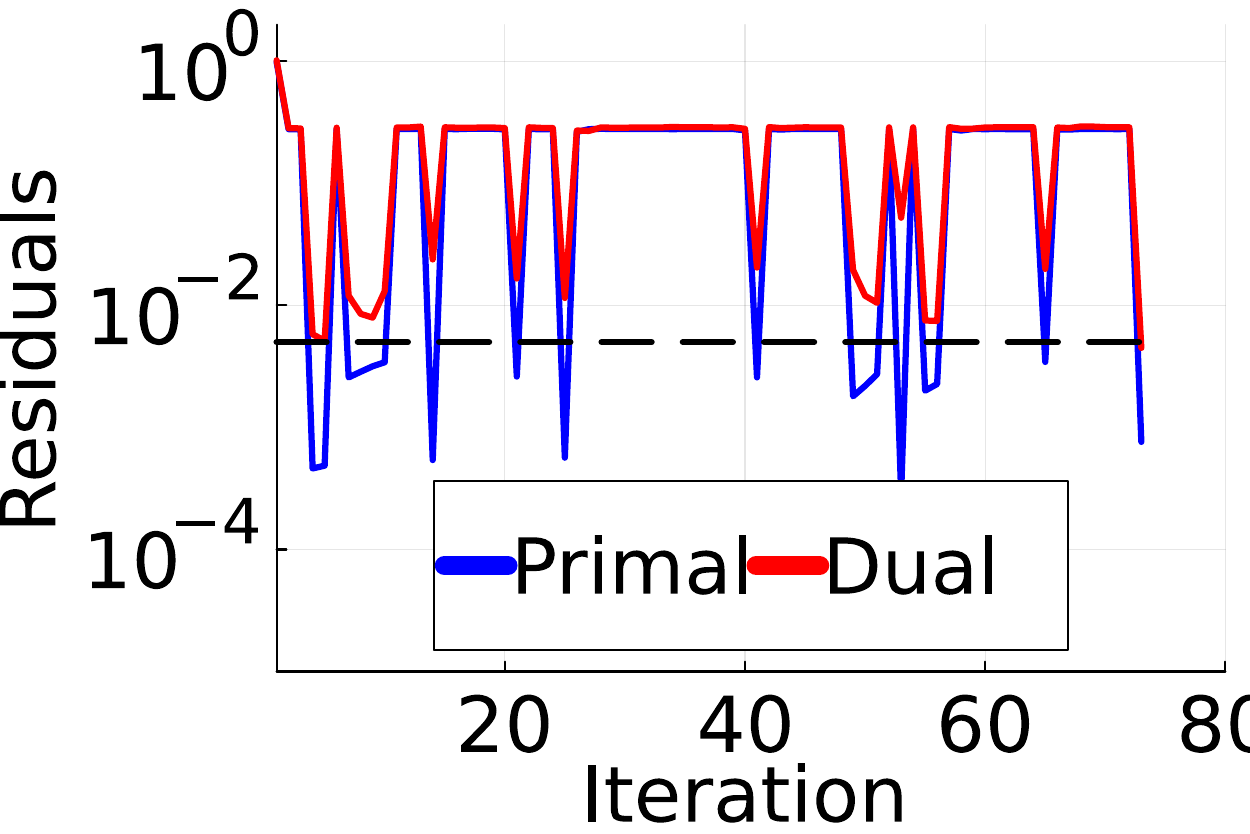}   
\includegraphics[width=.49\linewidth]{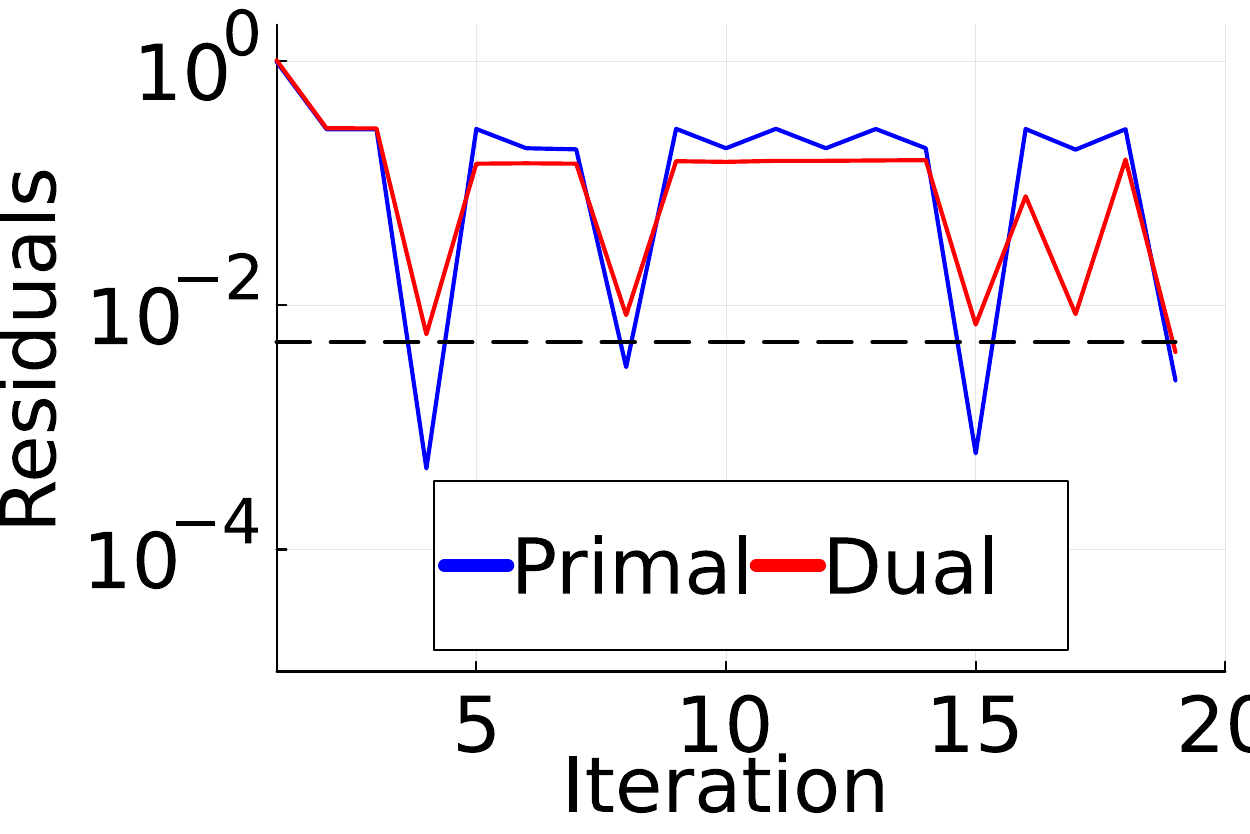} 
\caption{Effect of the NRB technique on the convergence.}
\label{fig:nrb_effect}
\end{figure}

To see the effect of the NRB technique on the convergence, we solved the UIUC-150 test instance when $E=5$ and report how the primal and dual residuals behave in Fig.~\ref{fig:nrb_effect}. 
Specifically, with a constant $\rho=10^2$, 73 iterations are consumed to solve the instance, while it takes only 19 iterations when NRB with $\rho^{(0)}=10^2$, $\beta=2$, and $\tau=10$ is used.

\subsection{Details on SL}

\begin{figure}[h!]
\centering      
\includegraphics[width=.49\linewidth]{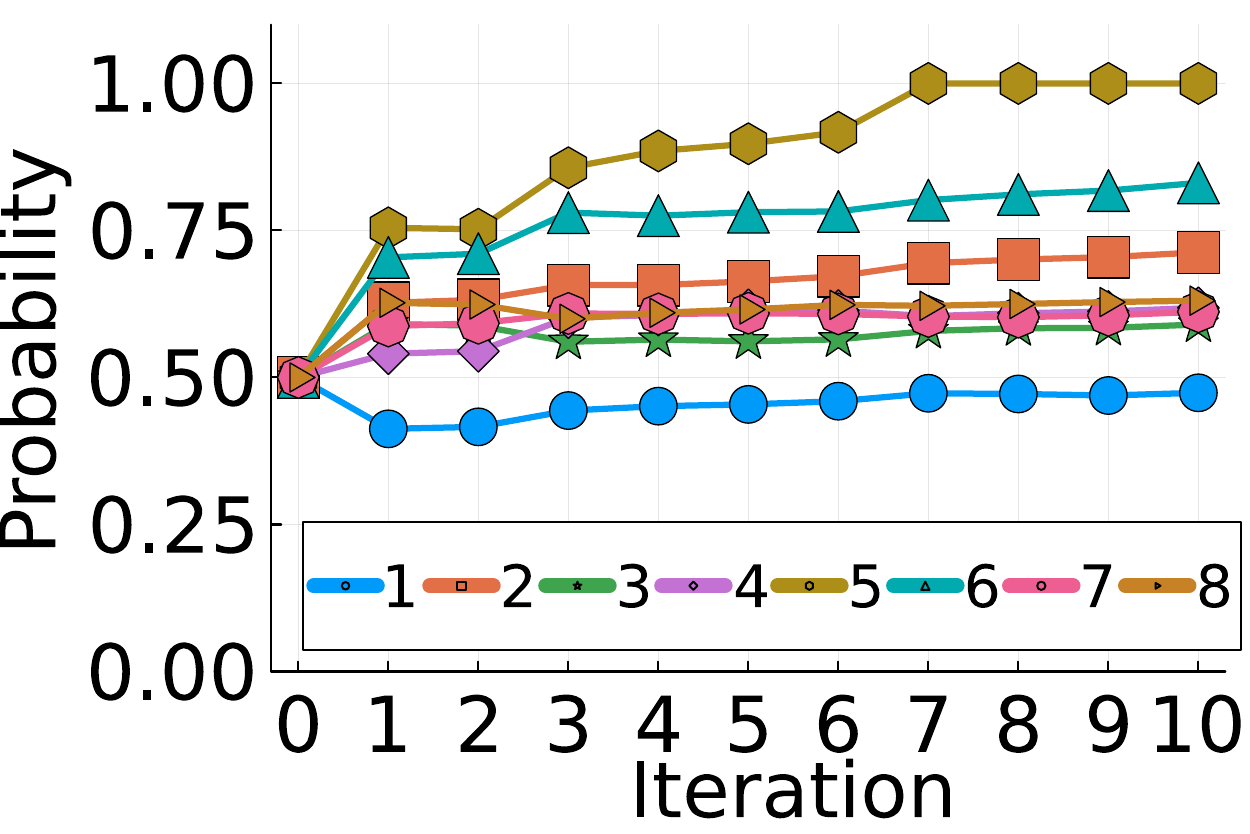} 
\includegraphics[width=.49\linewidth]{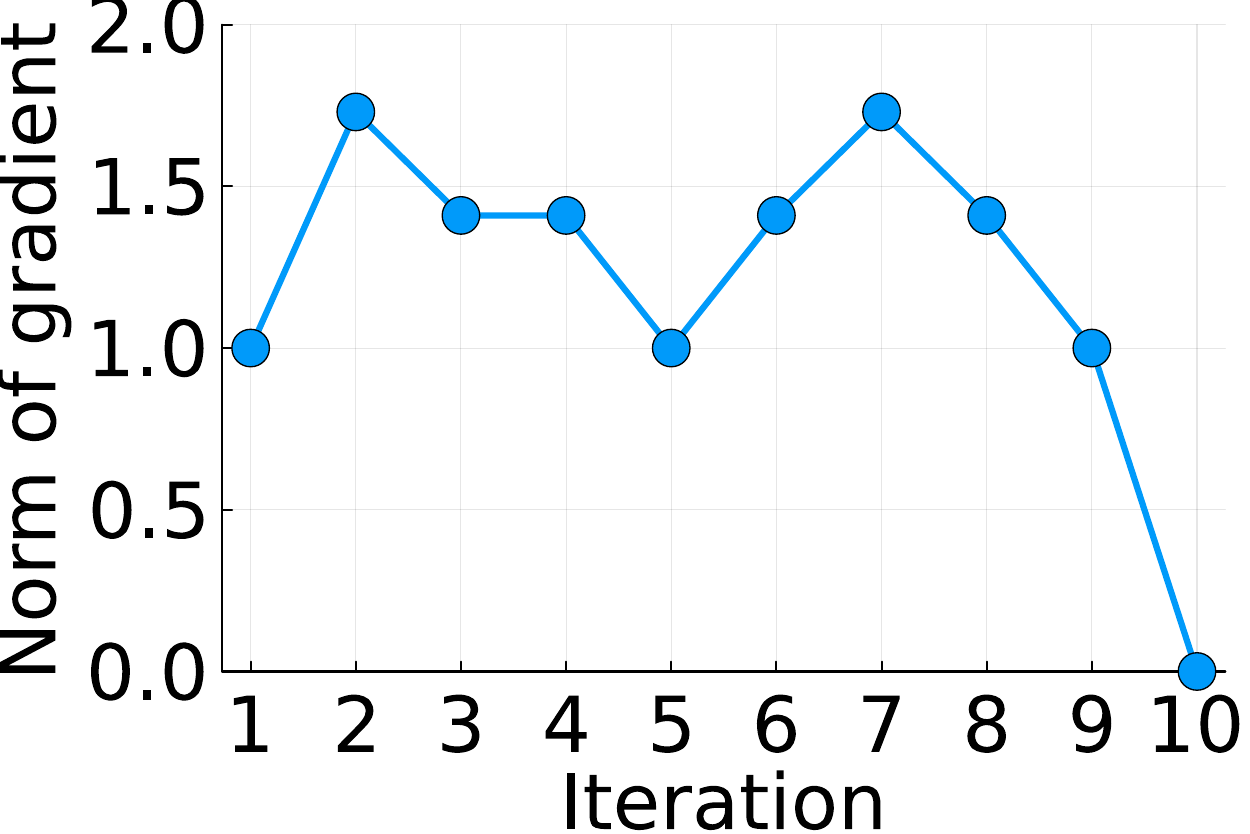} 
\caption{Variations of solutions (left) and the norm of gradients (right).}
\label{fig:sl_epri21}
\end{figure}

In this section we present details on the SL approach. 
First, we solved the EPRI-21 test instance, which has 8 substations, namely, candidate locations for installing blocking devices, when $E=5$ V/km. In Fig.~\ref{fig:sl_epri21} we report how the the solutions, namely, probabilities $\{p^{(t)}_1, \ldots, p^{(t)}_8\}$ in Eq.~\eqref{gradient_sl}, change over iterations and the norm of gradients (i.e., $\|g^{(t)}\|$ in Eq.~\eqref{gradient_sl}).
As reported, the solutions started from the initial values of $0.5$ change over iterations. For example, the probability of placing a blocking device at substation $8$ reaches $1$ upon termination, indicating that substation 8 should be chosen as an installation site.
For the remaining substations, the probabilities at the termination (i.e., 10th iteration where the norm of the gradient diminishes to zero) are used to determine whether or not to proceed with installation.

\begin{figure}[h!]
\centering    
  \includegraphics[width=.49\linewidth]{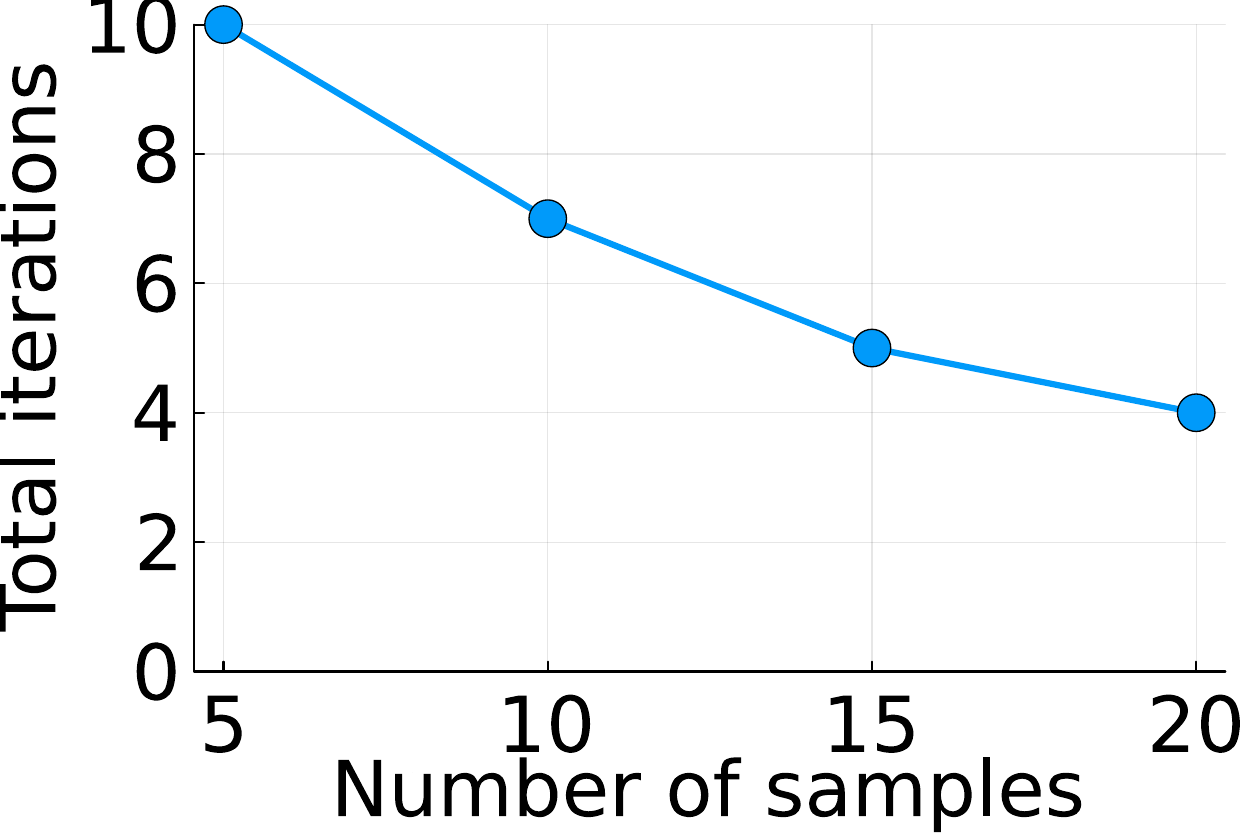}    
  \includegraphics[width=.49\linewidth]{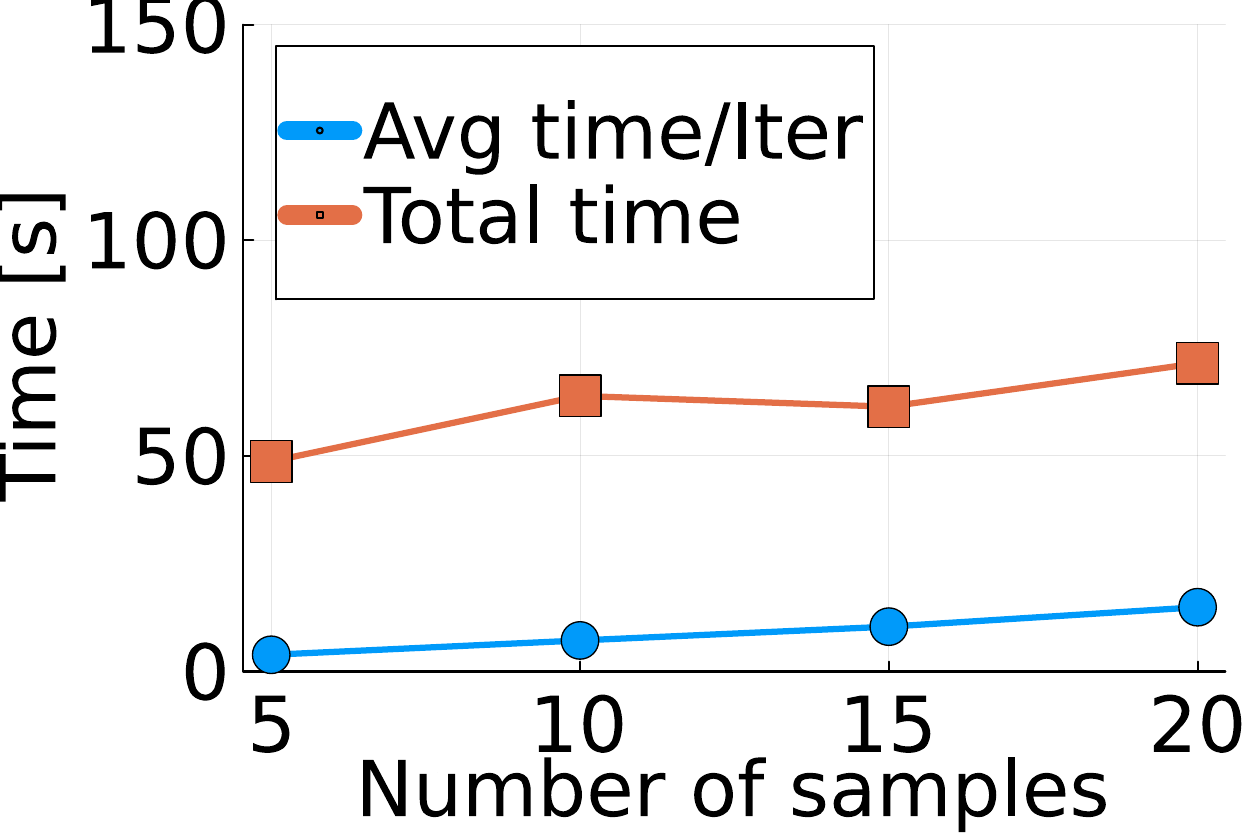}   
  \caption{Effect of the number $N$ of samples used in SL on computation.}
  \label{fig:sl_samples}
\end{figure}

The number, $N$, of samples can be linked to the batch size in the mini-batch SGD, a widely used learning algorithm. 
As $N$ increases, the algorithm's performance approaches that of a standard gradient descent algorithm.
With a larger $N$, however, the computation time per iteration increases since it entails computing the full or true gradient. 
In our case, where the objective is to obtain solutions within a 1-hour time limit, we restrict $N$ to be small.
In Fig.~\ref{fig:sl_samples} we illustrate the impact of $N$ on total iterations, average time per iteration, and total time. As anticipated, higher values of $N$ reduce the overall number of iterations needed for termination because of better solution quality per iteration, but this comes at the cost of longer computation time.
The choice of $N$ significantly influences the algorithm's performance, necessitating tuning based on the specific application's requirements.

\section{Conclusion and Future Work}

In this paper we derived a MINLP model that mathematically represents the GIC-BDP problem, and we proposed a heuristic algorithm by exploiting the structure of the problem.
As opposed to existing models that employ nonsmooth absolute-value equations to calculate effective GIC and intricate nonlinear, nonconvex equations involving trilinear terms with trigonometric functions to compute the amounts of power flow in transmission systems, we have formulated an alternative MINLP model to circumvent potential numerical instability concerns inherent in MINLP solvers.
This is achieved by employing a complementarity reformulation technique to smooth out the absolute-value equations.
Instead of solving the model using MINLP solvers, which often require substantial time due to the problem's inherent complexity, we have proposed a new heuristic algorithm, 3ADMM-B, that exploits the decomposition of the overall model into three segments. This approach renders each subproblem considerably simpler to solve. We also compared 3ADMM-B with a stochastic learning algorithm that optimizes the probability of GIC blocker placement. Compared with conventional solvers, the heuristics yield solutions of superior quality in significantly shorter timeframes. While 3ADMM-B is typically faster than the stochastic learning approach, the latter allows us to sample from the optimal probability distribution, allowing us to explore possible alternative solutions of similar quality.

As future work, we plan to integrate line switching decisions into the GIC-BDP problem, which serve to effectively counteract the detrimental effects of GIC on power grids by selectively deactivating specific transmission lines.
In contrast to the installation of blocking devices, line-switching decisions can be readily modified to enhance overall performance. 
To address this, we intend to structure the problem as a two-stage optimization process, where we determine the optimal placement of blocking devices in the first stage and make line-switching decisions in the second stage.
Given the prompt feasibility of line switching in practical scenarios, we intend to leverage machine learning methodologies to pretrain a model that will provide real-time recommendations on which lines to activate.
In consideration of this approach, the heuristic methods proposed in this paper can be employed to generate a collection of training datasets.

\section*{Acknowledgments}

The work was jointly funded by the U.S. Department of Energy's Office of Electricity Advanced Grid Modeling
program and the U.S. Department of Energy's Office of Science Scientific Discovery Through Advanced Computing (SciDAC) program under the project ``Space Weather Mitigation Planning.''
\bibliographystyle{unsrt}
\bibliography{reference.bib}

\begin{thebibliography}{10}

\bibitem{dhs}
\url{https://www.dhs.gov/science-and-technology/electromagnetic-pulse-empgeomagnetic-disturbance}.

\bibitem{NERC2012-gmd}
{NERC}.
\newblock {Effects of Geomagnetic Disturbances on the Bulk Power System System}.
\newblock Technical report, North American Electric Reliability Corporation, 2012.

\bibitem{barnes21-hiddenfailures}
A.~K.~{Barnes} et~al.
\newblock {The Risk of Hidden Failures to the United States Electrical Grid and Potential for Mitigation}.
\newblock In {\em {Proc. of the 2021 IEEE 53rd North American Power Symposium}}, Nov. 2021.

\bibitem{mate21-pmsgmd-cascade}
A.~{Mate} et~al.
\newblock {Relaxation Based Modeling of GMD Induced Cascading Failures in PowerModelsGMD.jl}.
\newblock In {\em {Proceedings of the 2021 IEEE/PES 53rd North American Power Symposium}}, 2021.

\bibitem{kappenman1991gic}
J.~G.~{Kappenman} et~al.
\newblock {GIC Mitigation: A Neutral Blocking/Bypass Device to Prevent the Flow of GIC in Power Systems}.
\newblock {\em IEEE Transactions on Power Delivery}, 6(3):1271--1281, 1991.

\bibitem{zhu2014blocking}
H.~{Zhu} et~al.
\newblock {Blocking Device Placement for Mitigating the Effects of Geomagnetically Induced Currents}.
\newblock {\em IEEE Trans. on Power Systems}, 30(4):2081--2089, 2014.

\bibitem{etemadi2014optimal}
A.~H.~{Etemadi} et~al.
\newblock {Optimal Placement of GIC Blocking Devices for Geomagnetic Disturbance Mitigation}.
\newblock {\em IEEE Transactions on Power Systems}, 29(6):2753--2762, 2014.

\bibitem{liang2015optimal}
Y.~{Liang} et~al.
\newblock {Optimal Blocker Placement for Mitigating the Effects of Geomagnetic Induced Currents Using Branch and Cut Algorithm}.
\newblock In {\em 2015 North American Power Symposium (NAPS)}, pages 1--6, 2015.

\bibitem{liang2019optimal}
Y.~{Liang} et~al.
\newblock {Optimal Blocking Device Placement for Geomagnetic Disturbance Mitigation}.
\newblock {\em IEEE Transactions on Power Delivery}, 34(6):2219--2231, 2019.

\bibitem{lu2017optimal}
M.~{Lu} et~al.
\newblock {Optimal Transmission Line Switching Under Geomagnetic Disturbances}.
\newblock {\em IEEE Trans. on Power Systems}, 33(3):2539--2550, 2017.

\bibitem{ryu2020algorithms}
M.~{Ryu} et~al.
\newblock {Algorithms for Mitigating the Effect of Uncertain Geomagnetic Disturbances in Electric Grids}.
\newblock {\em Electric Power Systems Research}, 189:106790, 2020.

\bibitem{ryu2022mitigating}
M.~{Ryu} et~al.
\newblock {Mitigating the Impacts of Uncertain Geomagnetic Disturbances on Electric Grids: A Distributionally Robust Optimization Approach}.
\newblock {\em IEEE Transactions on Power Systems}, 37(6):4258--4269, 2022.

\bibitem{ferris1997engineering}
M.~C.~{Ferris} et~al.
\newblock {Engineering and Economic Applications of Complementarity Problems}.
\newblock {\em Siam Review}, 39(4):669--713, 1997.

\bibitem{fletcher2004solving}
R.~{Fletcher} et~al.
\newblock {Solving Mathematical Programs with Complementarity Constraints as Nonlinear Programs}.
\newblock {\em Optimization Methods and Software}, 19(1):15--40, 2004.

\bibitem{scheel2000mathematical}
H.~{Scheel} et~al.
\newblock {Mathematical Programs with Complementarity Constraints: Stationarity, Optimality, and Sensitivity}.
\newblock {\em Mathematics of Operations Research}, 25(1):1--22, 2000.

\bibitem{raghunathan2005interior}
A.~U.~{Raghunathan} et~al.
\newblock {An Interior Point Method for Mathematical Programs with Complementarity Constraints (MPCCs)}.
\newblock {\em SIAM Journal on Optimization}, 15(3):720--750, 2005.

\bibitem{douglas1956numerical}
J.~{Douglas} et~al.
\newblock {On the Numerical Solution of Heat Conduction Problems in Two and Three Space Variables}.
\newblock {\em Transactions of the American mathematical Society}, 82(2):421--439, 1956.

\bibitem{eckstein1992douglas}
J.~{Eckstein} et~al.
\newblock {On the Douglas-—Rachford Splitting Method and the Proximal Point Algorithm for Maximal Monotone Operators}.
\newblock {\em Mathematical Programming}, 55:293--318, 1992.

\bibitem{wohlberg2017admm}
B.~{Wohlberg}.
\newblock {ADMM Penalty Parameter Selection by Residual Balancing}.
\newblock {\em arXiv preprint arXiv:1704.06209}, 2017.

\bibitem{attia2022stochastic}
A.~{Attia} et~al.
\newblock {Stochastic Learning Approach for Binary Optimization: Application to Bayesian Optimal Design of Experiments}.
\newblock {\em SIAM Journal on Scientific Computing}, 44(2):B395--B427, 2022.

\bibitem{BestuzhevaEtal2021ZR}
K.~{Bestuzheva} et~al.
\newblock {The SCIP Optimization Suite 8.0}.
\newblock ZIB-Report 21-41, Zuse Institute Berlin, December 2021.

\bibitem{kroger2018juniper}
O.~Kroger et~al.
\newblock {Juniper: An Open-Source Nonlinear Branch-and-Bound Solver in Julia}.
\newblock In {\em Integration of Constraint Programming, Artificial Intelligence, and Operations Research Conference}, pages 377--386, 2018.

\bibitem{Lubin2023}
M.~{Lubin} et~al.
\newblock {JuMP 1.0: Recent Improvements to a Modeling Language for Mathematical Optimization}.
\newblock {\em Mathematical Programming Computation}, 2023.

\end{thebibliography}

\iftrue
    \vspace{0.3cm}
    \begin{center}
        \scriptsize \framebox{
            \parbox{2.5in}{
                Government License (will be removed at publication):
                The submitted manuscript has been created by UChicago Argonne, LLC,
                Operator of Argonne National Laboratory (``Argonne").  Argonne, a
                U.S. Department of Energy Office of Science laboratory, is operated
                under Contract No. DE-AC02-06CH11357.  The U.S. Government retains for
                itself, and others acting on its behalf, a paid-up nonexclusive,
                irrevocable worldwide license in said article to reproduce, prepare
                derivative works, distribute copies to the public, and perform
                publicly and display publicly, by or on behalf of the Government. The Department of Energy will provide public access to these results of federally sponsored research in accordance with the DOE Public Access Plan. http://energy.gov/downloads/doe-public-access-plan.
            }
        }
        \normalsize
    \end{center}
\fi

\end{document}